\documentclass[11pt]{article}
\usepackage{amssymb,amsfonts,amsmath,latexsym,epsf,tikz,url}

\newtheorem{theorem}{Theorem}[section]

\newtheorem{corollary}[theorem]{Corollary}
\newtheorem{lemma}[theorem]{Lemma}

\newcommand{\proof}{\noindent{\bf Proof.\ }}
\newcommand{\qed}{\hfill $\square$\medskip}

\begin{document}

\title{Counting the number of weakly connected dominating sets of graphs}

\author{
Saeid Alikhani$^{}$\footnote{Corresponding author}, Somayeh Jahari, Mohammad Mehryar
}

\date{\today}

\maketitle

\begin{center}
Department of Mathematics, Yazd University, 89195-741, Yazd, Iran\\
{\tt alikhani@yazd.ac.ir}
\end{center}


\begin{abstract}
Let $G=(V(G),E(G))$ be a simple graph. A non-empty set $S\subseteq V (G)$ is a weakly connected dominating set in $G$, if
the subgraph obtained from $G$ by removing all edges each joining any two vertices in $V (G)\setminus S$ is
connected.
In this paper, we consider some graphs and study the number of their weakly connected dominating sets.

\end{abstract}


\medskip
\noindent{\bf AMS Subj.\ Class.}:  05C05, 05C69

\section{Introduction} \label{introduction}

Let $G=(V(G),E(G))$ be a simple connected graph. A set $S\subseteq V(G)$ is a {\it dominating set} if every vertex in $V(G)\backslash S$ is adjacent to at least one vertex in $S$.
The {\it domination number} $\gamma(G)$ is the minimum cardinality of a dominating set in $G$. There are various domination numbers in the literature (\cite{dom}).
For a detailed treatment of domination theory, the reader is referred to \cite{domination}. A non-empty $S\subseteq V (G)$, $S$ is called
a weakly connected dominating set (w.c.d.s.) of $G$, if the subgraph obtained from
$G$ by removing all edges each joining any two vertices in $V (G)\setminus S$ is
connected. The weakly connected domination number $\gamma_w(G)$, is defined
to be the minimum integer $k$ with $|S| = k$ for some weakly connected
dominating set $S$ of $G$ (see \cite{dom,koh}).

A dominating set
with cardinality $\gamma_w(G)$ is  called a {\it $\gamma_w$-set}.    Let ${\cal D}_w(G,i)$ be the family of
weakly connected  dominating sets of a graph $G$ with cardinality $i$ and let
$d_w(G,i)=|{\cal D}_w(G,i)|$.
The number of dominating sets of a graph has been actively studied in recent
years (\cite{euro,utilitas,ars,elec}). In this paper, we shall count the number of weakly connected dominating sets of a graph $G$.

For two graphs $G = (V,E)$ and $H=(W,F)$, the corona $G\circ H$ is the graph arising from the
disjoint union of $G$ with $| V |$ copies of $H$, by adding edges between
the $i$th vertex of $G$ and all vertices of $i$th copy of $H$ (\cite{Fruc}).  The join $ G + H$ of two graph $G$ and $H$ with disjoint vertex sets $V$ and $W$ and
edge sets $E$ and $F$ is the graph union $G\cup H$ together with all the edges joining $V$ and
$W$.
\medskip

In the next section, we consider specific graphs and
study the number of their weakly connected dominating sets. In Section $3$, we consider graphs with specific
construction, denoted by $G(m)$ and
construct all their weakly connected dominating sets.  As an example of these graphs,  we study the structure of weakly connected dominating sets and the number of weakly connected dominating sets of paths. Finally, we study the number of weakly connected dominating sets of cycles in the last
section.

As usual, we use $\lceil x \rceil$, $\lfloor x\rfloor$ for the
smallest integer greater than or equal to $x$ and the largest
integer less than or equal to $x$, respectively. Also we denote the complete graph, path and cycle of order $n$ by $K_n$, $P_n$ and $C_n$, respectively. Also $K_{1,n}$ is the star graph with $n+1$ vertices. 
In this article,
we denote $\{1,2,\ldots,n\}$ simply by $[n]$.

\section{Weakly connected dominating  sets of specific graphs}

In this section we consider specific graphs and study their weakly connected dominating sets with cardinality $i$, for $\gamma_w(G)\leq i \leq |V(G)|$. It is
well-known and generally accepted that the problem of determining
the domination number and dominating sets (and so weakly connected domination number and weakly connected dominating sets) of an arbitrary graph is
difficult. Since this problem has been shown to be
NP-complete (see~\cite{johnson}),
we shall consider  in this  section, specific graphs.

First we consider the complete graph $K_n$ and the star graphs $K_{1,n}$.
The number of weakly connected dominating sets of $K_n$ and $K_{1,n}$ are easy to compute.

\begin{theorem}
	\begin{enumerate}
		\item[(i)] For every $n\in \mathbb{N}$, and $1\leq i\leq n$, $d_w(K_n,i)={n\choose i}$.
		\item[(ii)] For every $1\leq i \leq n-1$, $d_w(K_{1,n},i)={n\choose i-1}$.
		\item[(iii)] $d_w(K_{1,n},n)=n+1$ and $d_w(K_{1,n},n+1)=1$.
	\end{enumerate}
\end{theorem}

The following theorems give the weakly connected domination number of corona and join of  two graphs:

\begin{theorem} \rm\cite{IMF}
	Let $G$ be a connected graph with $|V(G)|\geq 2$ and $H$ an arbitrary  graph.
	Then $\gamma_w(G\circ H) = |V (G)|$.
\end{theorem}

\begin{theorem}\rm\cite{IMF}
	For two graphs $G$ and $H$,
	$$\gamma_w(G+H)=\left\{
	\begin{array}{lcl}
	1  & \mbox{ \rm if $\gamma(G)=1$ or $\gamma(H)=1$}; \\[10pt]
	2 & \mbox{ \rm otherwise}.
	\end{array}
	\right.$$
\end{theorem}

The following theorem gives the number of w.c.d.s. of $G_1+G_2$.
\begin{theorem}\label{join}
	Let $G_1$ and $G_2$ be  connected graphs of order  $n_1$ and $n_2$,
	respectively. Then, for two natural numbers $i_1, i_2$, and $i\geq \gamma_w(G_1+G_2)$,
	$$d_w(G_1+G_2,i)=d_w(G_1,i)+d_w(G_2,i)+\displaystyle\sum_{i_1+i_2=i} {n_1\choose i_1}{n_2\choose i_2}.$$
\end{theorem}
\proof.  Let $i$ be a natural number  $1\leq i \leq
n_1+n_2$.  We want to determine $d_w(G_1+G_2,i)$.  If $i_1$ and
$i_2$ are two natural numbers such that $i_1+i_2=i$ ,
then clearly,  for every $D_1\subseteq  V(G_1)$ and $D_2\subseteq V(G_2)$ ,
such that $|D_j|=i_j$ ,  $j=1,2$ ,  $D_1\cup D_2$ is a weakly connected dominating set
of $G_1+ G_2$.  Moreover, if  $D\in  {\cal D}_w(G_1,i)$,
then $D$ is a weakly connected dominating set for $G_1+G_2$ of size $i$.  The
same is true for every $D\in {\cal D}_w(G_2,i)$. Therefore
we have the result. \qed

The following corollary gives the  relationship between the number of w.c.d.s.  of wheels $W_n$ and cycles $C_n$:

\begin{corollary}
	For every $n\geq 4$, $d_w(W_n,i)= \left\{
	\begin{array}{lcl}
	1  & \mbox{ \rm if $i=1$}; \\[10pt]
	d_w(C_{n-1},i)+{n-1\choose i-1}  & \mbox{ \rm if $i\geq 2$}.
	\end{array}
	\right.
	$
\end{corollary}
\proof
Since $W_n=C_{n-1}+K_1$, by Theorem \ref{join} we have,
$$d_w(W_n,i)= \left\{
\begin{array}{lcl}
d_w(C_{n-1},i)+1  & \mbox{ \rm if $i=1$}; \\[10pt]
d_w(C_{n-1},i)+{n-1\choose i-1}  & \mbox{ \rm otherwise}.
\end{array}
\right.
$$
Since for every $n\geq 4$, $d_w(C_{n-1},1)=0$, we have the result.\qed

\section{Weakly connected dominating sets of $G(m)$}

In this section, we shall study the weakly connected dominating sets (w.c.d.s.) of specific graphs denoted by $G(m)$. As an example of graphs $G(m)$, we construct w.c.d.s. of paths and count the number of w.c.d.s. of paths.

\begin{figure}[h]
	\hspace{3cm}
	\includegraphics[width=5.7cm,height=0.6cm]{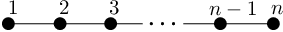}
	\caption{\label{figure3} The path $P_n$ with vertices labeled
		$[n]$.}
\end{figure}

A path is a connected graph in which two vertices have degree
one and the remaining vertices have degree two. Let $P_n$ be the
path with $V(P_n)=[n]$ and
$E(P_n)=\{\{1,2\},\{2,3\},\ldots,\{n-1,n\}\}$, see Figure~\ref{figure3}.


Let $P_{m+1}$ be a path with vertices labeled by $y_0, y_1,
\ldots, y_m$, for $m\geq 0$ and let $v_0$ be a specific vertex of
a graph $G$. Denote by $G_{v_0} (m)$ a graph obtained from $G$ by
identifying the vertex $v_0$ of $G$ with an end vertex $y_0$ of
$P_{m+1}$. It is clear that if the path is glued to a different
vertex $v_1$ of $G$, then the two graphs $G_{v_1} (m)$ and
$G_{v_0} (m)$ may not be isomorphic. It depends on the vertex to
which  we glue the path. If throughout our discussion, this vertex
is fixed, then we shall simply use the notation $G(m)$ (if there
is no likelihood of confusion).

We need the following  lemma  to obtain our main results in this section:

\begin{lemma}\label{lemma1}
	\begin{enumerate}
		\item[\rm (i)] ${\cal D}_w(G(m),i)=\emptyset$ if and only if
		$i>|V(G(m))|$ or $i<\gamma_w(G(m))$,
		\item[\rm (ii)] If $e\in E(G)$,
		then $\gamma_w(G-e)-1\leq \gamma_w(G) \leq \gamma_w(G-e)$, {\rm
			(see~\cite{czech})}
		\item[\rm(iii)] For any $m\in \mathbb{N}$,
		$\gamma_w(G(m-1))\leq \gamma_w(G(m))\leq \gamma_w(G(m-1))+1$. \rm(by
		$(ii)$ above).
	\end{enumerate}
\end{lemma}

We need the following easy lemma  and theorem:

\begin{lemma}\label{lemma3}
	For every $n\in \mathbb{N}$,
	$\gamma_w(P_{n})=\lfloor\frac{n}{2}\rfloor$.
\end{lemma}

\begin{theorem}\label{theorem3.4.1}
	Suppose that $G(m)$ is the graph defined in this section.
	Then for every $m\geq 0$,
	\[
	\gamma_w(G(m))=\left\{
	\begin{array}{lcl}
	\gamma_w(G)+\lfloor\frac{m-1}{2}\rfloor  & \mbox{ \rm if $G$ has a $\gamma$-set containing $y_0$}; \\[10pt]
	\gamma_w(G)+\lfloor\frac{m}{2}\rfloor  & \mbox{ \rm otherwise}.
	\end{array}
	\right.
	\]
\end{theorem}
\proof  If $y_0$ is in the $\gamma_w$-set of $G$, then to
obtain the $\gamma_w$-set
of $G(m)$ it suffices to dominate the path with vertices $\{y_2,\ldots,y_m\}$,
otherwise we dominate the path  with vertices
$\{y_1,y_2,\ldots,y_m\}$. Therefore;  by Lemma~\ref{lemma3}, the
proof is complete.\qed

To enumerate  the  weakly  connected dominating set of $G(m)$ with cardinality
$i$, no need to consider w.c.d.s. of $G(m - 3)$ with
cardinality $i - 1$. 
Therefore, we only need to consider w.c.d.s.  in $G(m - 1)$ and
$G(m - 2)$ with cardinality $i-1$. The families of
these weakly connected dominating sets (w.c.d.s.) can be empty or otherwise. Thus, we have
four cases of whether these two families are empty or not. We do not need to consider the case that
${\cal D}_w(G(m - 1), i - 1) ={\cal  D}_w(G(m - 2), i - 1) = \emptyset$,
because it implies ${\cal D}_w(G(m), i) =\emptyset$. Also the case ${\cal D}_w(G(m - 1), i - 1)\neq \emptyset$, ${\cal  D}_w(G(m - 2), i - 1) = \emptyset$ does not exists. Thus, we only need to
consider two cases. We consider these cases in
Theorem~\ref{theorem2} which construct the w.c.d.s.   of $G(m)$.

\begin{theorem}\label{theorem2}
	\begin{enumerate} \item[\rm(i)] If  ${\cal D}_w(G(m-1),i-1)=\emptyset$ and ${\cal D}_w(G(m-2),i-1)\neq
		\emptyset$,  then ${\cal D}_w(G(m),i)=\Big\{\{y_{m-1}\}\cup
		X\,|\,X\in{\cal D}_w(G(m-2),i-1)\Big\}$,
		
		\item[\rm(ii)] If  ${\cal
			D}_w(G(m-2),i-1)\neq \emptyset$, ${\cal D}_w(G(m-1),i-1)\neq\emptyset$, then ${\cal D}_w(G(m),i)=\Big\{\,\{y_m\}\cup X_1, \{y_{m-1}\}\cup
		X_2\,|\,X_1\in
		{\cal D}_w(G(m-1),i-1),X_2\in {\cal D}_w(G(m-2),i-1) \Big\}$
		
	\end{enumerate}
\end{theorem}

\proof 
\begin{enumerate}
	\item[(i)] Obviously  $\Big\{\{y_{m-1}\}\cup
	X\,|\,X\in{\cal D}_w(G(m-2),i-1)\Big\}\subseteq {\cal D}_w(G(m),i)$. Now suppose that
	$Y\in {\cal D}_w(G(m),i)$. Then at least one of the vertices $y_m$
	or $y_{m-1}$  is in $Y$. If $y_m\in Y$ then at least one of the vertices $y_{m-1}$ or $y_{m-2}$ is in
	$Y$. If $y_{m-1}\in Y$, then $Y-\{y_m\}\in {\cal
		D}_w(G(m-1),i-1)$ a contradiction. So $y_{m-2}\in Y$ and $Y-\{y_{m-1}\}\in {\cal D}_w(G(m-2),i-1)$. Therefore ${\cal D}_w(G(m),i)\subseteq \Big\{\{y_{m-1}\}\cup
	X\,|\,X\in{\cal D}_w(G(m-2),i-1)\Big\}$.
	
	\item[(ii)]
	Obviously $\Big\{\,\{y_m\}\cup X_1, \{y_{m-1}\}\cup
	X_2\,|\,X_1\in
	{\cal D}_w(G(m-1),i-1),X_2\in {\cal D}_w(G(m-2),i-1) \Big\}\subseteq {\cal D}_w(G(m),i)$.

	Now, let $Y \in {\cal D}_w(G(m),i)$, then $y_m \in Y$ or $y_{m-1}
	\in Y$. If $y_m \in Y$, then  at least one
	vertex labeled $y_{m-1}$ or $y_{m-2}$ is in $Y$. If
	$y_{m-1}\in Y$, then $Y=X\cup \{y_m\}$ for some
	$X\in {\cal D}(G(m-1),i-1)$. If $y_{m-2}\in Y$,
	then $Y=X\cup \{y_{m-1}\}$ for some $X\in
	{\cal D}(G(m-2),i-1)$. So we have the result.\qed

\end{enumerate}

\begin{theorem}\label{theorem3'}
	For every $m\geq 2$,
	\[
	d_w(G(m),i)=d_w(G(m-1),i-1)+d_w(G(m-2),i-1).
	\]
\end{theorem}

\proof  It follows from  Theorem~\ref{theorem2}.\qed

Since $P_n=P_1(n-1)$, we can apply the results for
the graph $G(m)$  to obtain some properties
of w.c.d.s.  and their numbers  for  paths. We denote
${\cal D}_w(P_{n},i)$ simply by ${\cal P}_{n}^{i}$.

For the construction of ${\cal P}_n^i$, by Theorem~\ref{theorem2}, we only need to consider two  families
${\cal P}_{n-1}^{i-1}$ and ${\cal P}_{n-2}^{i-1}$.

\begin{theorem}\label{theorem5}
	For every $n\geq 3$ and $i\geq \lfloor\frac{n}{2}\rfloor$,
	\begin{enumerate}
		\item[(i)] If ${\cal P}_{n-1}^{i-1}=\emptyset$ and ${\cal P}_{n-2}^{i-1}\neq\emptyset$,
		then ${\cal P}_n^i=\Big\{ X\cup\{n-1\}|\;X\in{\cal P}_{n-2}^{i-1}\Big\}$.
		
		\item[(ii)]  If  ${\cal P}_{n-1}^{i-1}\neq\emptyset$ and $ {\cal P}_{n-2}^{i-1}\neq\emptyset$, then
		$${\cal P}_{n}^i=\Big\{\{n\}\cup X_1, \{n-1\}\cup X_2\,|\,X_1\in
		{\cal P}_{n-1}^{i-1},X_2\in {\cal P}_{n-2}^{i-1} \Big\}.$$
		
	\end{enumerate}
\end{theorem}

\proof It follows from Theorem \ref{theorem2}.\qed


The following theorem gives a recurrence for the number of w.c.d.s. of $P_n$.
\begin{theorem}\label{theorem3"}
	For every $n\geq 3$ and $\lfloor\frac{n}{2}\rfloor \leq i \leq n$, $d_w(P_n,i)=d_w(P_{n-1},i-1)+d_w(P_{n-2},i-1)$, with initial values
	$d_w(P_1,1)=1$, $d_w(P_2,1)=2$ and $d_w(P_2,2)=1$.
\end{theorem}
\proof It follows from Theorem \ref{theorem3'}.\qed

Using  Theorem~\ref{theorem3"} ,  we obtain  $d_w(P_{n},j)=|{\cal P}_n^j|$ for $1\leq n\leq 10$ in Table 1.

\begin{center}
	\begin{footnotesize}
		\small
		\begin{tabular}{r||ccccccccccc}
			$j$&$1$&$2$&$3$&$4$&$5$&$6$&$7$&$8$&$9$&$10$&\\ [0.1ex]
			\hline
			$d_w(P_1,j)$&1&&&&&&&&&&\\
			$d_w(P_2,j)$&2&1&&&&&&&&&\\
			$d_w(P_3,j)$&1&3&1&&&&&&&&\\
			$d_w(P_4,j)$&&3&4&1&&&&&&&\\
			$d_w(P_5,j)$&&1&6&5&1&&&&&&\\
			$d_w(P_6,j)$&&&4&10&6&1&&&&&\\
			$d_w(P_7,j)$&&&1&10&15&7&1&&&&\\
			$d_w(P_8,j)$&&&&5&20&21&8&1&&&\\
			$d_w(P_9,j)$&&&&1&15&35&28&9&1&&\\
			$d_w(P_{10},j)$&&&&&6&35&56&36&10&1&
		\end{tabular}
	\end{footnotesize}
\end{center}
\begin{center}
	{Table 1.} $d_w(P_{n},j)$, the number of w.c.d.s.  of $P_n$ with cardinality $j$.
\end{center}

Here, we shall  solve the recurrence relation with two variables for $d_w(P_n,j)$ in Theorem \ref{theorem3"}. Corresponding to this recurrence relation, we state an elementary combinatorial problem.

Suppose that we have $n$ boxes in the row and $j$ objects. We want to count the number of
permutations of these items in boxes such that there is at most one object  in each box and no two adjacent boxes can be empty. It is easy to see that the
answer of this problem  is ${j+1 \choose n-j}$. We can see that if
$a_{n,j}$ is the solution of this problem, then we have the following recurrence relation with these initial values $a_{1,1}=1$, $a_{2,1}=2$ and $a_{2,2}=1$:
\[
a_{n,j}=a_{n-1,j-1}+a_{n-2,j-1}.
\]
So we have the following result:

\begin{theorem}
	For every $n\in \mathbb{N}$ and $\lfloor\frac{n}{2}\rfloor\leq j\leq n$, $d_w(P_n,j)={j+1\choose n-j}$.
\end{theorem}

\section{Weakly connected dominating sets of $C_n$}

\begin{figure}[h]
	\hspace{4.3cm}
	\includegraphics[width=3.2cm, height=3.cm]{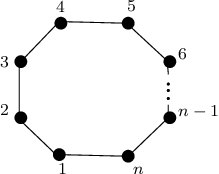}
	\caption{\label{figure2'} The cycle $C_n$ with vertices labeled $[n]$.}
\end{figure}

Let $C_n, n\geq 3$, be the cycle with $n$ vertices $V(C_n)=[n]$ and
$E(C_n)=\{\{1,2\},\{2,3\},\ldots,\{n-1,n\},\{n,1\}\}$, see Figure~\ref{figure2'}.
In this section, we consider the number of weakly dominating sets of cycle $C_n$.  Using  Maple programme,  we obtain  $d_w(C_{n},j)=|{\cal C}_n^j|$ for $1\leq n\leq 14$ in Table 2.

\begin{center}
	\begin{footnotesize}
		\small
		\begin{tabular}{r||cccccccccccccc}
			$j$&$1$&$2$&$3$&$4$&$5$&$6$&$7$&$8$&$9$&$10$&$11$&12&13&14\\ [0.1ex]
			\hline
			$d_w(C_1,j)$&1&& & & & & & & & & & &&\\
			$d_w(C_2,j)$&2&1& & & & & & & & & && &\\
			$d_w(C_3,j)$&3&3&1 & & & & & & & & & &&\\
			$d_w(C_4,j)$&&6&4 & 1& & & & & & & & &&\\
			$d_w(C_5,j)$&&5&10&5 &1 & & & & & & & &&\\
			$d_w(C_6,j)$& && 14&15 &6 &1 & & & & & & &&\\
			$d_w(C_7,j)$& &&7 &28 &21 &7 &1 & & & & & &&\\
			$d_w(C_8,j)$& && & 26& 48&28 &8 &1 & & & & &&\\
			$d_w(C_9,j)$& && & 9& 63& 75&36 &9 &1 & & & &&\\
			$d_w(C_{10},j)$& && & &42 & 125& 110&45 &10 &1 & & &&\\
			$d_w(C_{11},j)$& && & & 11& 121& 220&154 & 55&11 &1 & &&\\
			$d_w(C_{12},j)$& && & & &62 &276 &357 &208 &66 & 12&1 &&\\
			$d_w(C_{13},j)$& && & & & 13&208 &546 &546 & 273&78 & 13&1&\\
			$d_w(C_{14},j)$& && & & && 86&539 &980 &798 & 350&91 & 14&1\\
		\end{tabular}
	\end{footnotesize}
\end{center}
\begin{center}
	{Table 2.} $d_w(C_{n},j)$, the number of weakly connected dominating sets of $C_n$ with cardinality $j$.
\end{center}

\begin{lemma}\label{lemma2}
	For every $n\in \mathbb{N}$,
	$\gamma_w(C_{n})=\lfloor\frac{n}{2}\rfloor$.
\end{lemma}
\proof
We consider graph $C_n$ with the vertex set $V(C_n)=\{v_1, v_{2}, \dots, v_{n}\}$ and the edge set
$E(C_n)=\{\{v_{i-1}, v_{i}\}~|~1\leq i \leq n~and~ v_0=v_n\}$.
we have $\gamma_w(C_n) = 1$ for $n\leq 3$.  Now assume $n \geq 4$ and use induction on $n$. Suppose  that $S$ is a minimum weakly dominating
set of $C_n$ and consider the vertex  $v_i \in S$.  Since $N(v_i)=\{v_{i-1},v_{i+1}\}$ and
$C_n - \{v_i, v_{i-1}\}$ is a path with $n - 2$ vertices, by induction, $\gamma_w(C_n) = 1 + \gamma_w(P_{n-2}) = 1 + \lfloor \frac{n-2}{2} \rfloor = \lfloor \frac{n}{2} \rfloor$. \qed

The following theorem gives the number of w.c.d.s. of $C_n$ with cardinality $n-3\leq i \leq n$.
\begin{lemma}\label{lem3}
	\begin{enumerate}
		\item[(i)]
		For every $n\geq 4$, and $n-2 \leq i \leq n$, $d_w(C_n,i)={n\choose i}$.
		\item[(ii)]
		For every $n\geq 6$,  $d_w(C_n,n-3)=\frac{(n+1)n(n-4)}{6}$.
	\end{enumerate}
\end{lemma}
\proof
\begin{enumerate}
	\item[(i)]
	Since  $d_w(C_n,n)=1$ and $d_w(C_n,n-1)=n$, so the result is true for $i\in\{n-1,n\}$. We prove $d_w(C_n,n-2)={n\choose n-2}$.
	Clearly   for every $S\in \mathcal{D}_w(C_{n-1},n-2)$ we have $S\in \mathcal{D}_w(C_n,n-2)$.
	Also if $S_1 \in \mathcal{D}_w(C_{n-1},n-3)$, then $S_1\cup \{n\}\in \mathcal{D}_w(C_n,n-2)$.
	On the other hand, if $S\in \mathcal{D}_w(C_n,n-2)$, and $n \in S$, then $S\setminus \{n\} \in \mathcal{D}_w(C_{n-1},n-3)$. If $n\not\in S$, then $S\in \mathcal{D}_w(C_{n-1},n-2)$.
	Therefore
	\[ d_w(C_n,n-2)=d_w(C_{n-1},n-3)+d_w(C_{n-1},n-2).\]
	Now using induction, we have,
	\begin{eqnarray*}
		d_w(C_n,n-2)={n \choose n-2}.
	\end{eqnarray*}

	\item[(ii)] First we prove  that,
	\begin{eqnarray}
	d_w(C_n,n-3)=d_w(C_{n-1},n-4)+d_w(C_{n-1},n-3)-1.
	\end{eqnarray}
	Clearly every  w.c.d.s.  $S$ of  $C_n$ with  cardinality  $n-3$ is a w.c.d.s.  of $C_{n-1}$ with cardinality  $n-3$, except for the following cases:
	\begin{itemize}
		\item the set $S$  contains vertices $1$ and $ 2$,
		\item  the set $S$ contains vertices $1$ and $ n-1$,
		\item  the set $S$  contains vertices $n-1$ and $ n-2$.
	\end{itemize}
	Also it is easy to see that if $S\in \mathcal{D}_w(C_{n-1},n-4)$ or $S$ is any of the sets $\{2, \cdots, n-3\}$ and $\{3, \cdots, n-2\}$, then
	$S\cup \{n\}$ is a  w.c.d.s.  of size $n-3$ in cycle $C_{n}$.
	Consequently
	\begin{eqnarray*}
		d_w(C_n,n-3)&=&d_w(C_{n-1},n-4)+2+d_w(C_{n-1},n-3)-3 \\
		&=&d_w(C_{n-1},n-4)+d_w(C_{n-1},n-3)-1.
	\end{eqnarray*}
	Using equation (1),  we have,
	$$d_w(C_n,n-3)=\sum_{i=5}^n \frac{(i-1)(i-2)}{2}-(n-4), $$ and by easy computation we have the result.
	\quad\qed\
\end{enumerate}

\end{document}